\documentclass{article}

\usepackage[french,british]{babel}

\usepackage[T1]{fontenc}

\usepackage{amsmath}

\usepackage{amsfonts}
\usepackage{amsbsy}
\usepackage{amssymb}

\usepackage{amsthm}

\theoremstyle{plain}

\usepackage{latexsym}
\usepackage{amsmath}
\usepackage{amsfonts}
\usepackage{amssymb}
\usepackage{psfrag}
\usepackage{graphicx}

\usepackage{amsmath,amssymb}
\usepackage{graphicx}

\usepackage{amsmath,amssymb}
\usepackage{graphicx}

\newtheorem{ozn}{Definition}[section]
\newtheorem{thm}{Theorem}[section]

\newtheorem{nas}{Corollary}[section]

\newtheorem{lema}{Lemma}[section]

\newcommand{\eps}{\varepsilon}
\newcommand{\me}{\mathbf}
\newcommand{\mr}{\mathbb}
\newcommand{\mt}{\mathsf}
\newcommand{\md}{\mathcal}
\newcommand{\ld}{\left}
\newcommand{\rd}{\right}

\newcommand{\ip}{\int_{-\pi}^{\pi}}

\newcommand{\be}{\begin{equation}}
\newcommand{\ee}{\end{equation}}

\newcommand{\bem}{\begin{multline}}
\newcommand{\eem}{\end{multline}}

\newcommand{\bml}{\begin{multline*}}
\newcommand{\eml}{\end{multline*}}

\newcommand{\beg}{\begin{gather}}
\newcommand{\eeg}{\end{gather}}

\begin{document}

\title{Minimax-robust estimation problems for sequences with periodically stationary increments observed with noise}

\author{
Maksym Luz\thanks {BNP Paribas Cardif, Kyiv, Ukraine, maksym.luz@gmail.com},
Mikhail Moklyachuk\thanks
{Department of Probability Theory, Statistics and Actuarial
Mathematics, Taras Shevchenko National University of Kyiv, Kyiv 01601, Ukraine, moklyachuk@gmail.com}
}

\date{\today}

\maketitle

\renewcommand{\abstractname}{Abstract}
\begin{abstract}
 The problem of optimal estimation of linear functionals
constructed from the unobserved values of a stochastic sequence with periodically stationary increments based on observations of the sequence with stationary noise is considered.
For sequences with known spectral densities, we obtain formulas for calculating values of the mean square errors and the spectral characteristics of the optimal estimates of the functionals.
Formulas that determine the least favorable spectral densities and the minimax-robust spectral
characteristics of the optimal linear estimates of functionals
are proposed in the case where spectral densities of the sequence
are not exactly known while some sets of admissible spectral densities are specified.
\end{abstract}

\vspace{2ex}
\textbf{Keywords}:{periodically stationary increments, minimax-robust estimate, least favorable spectral density, minimax-robust spectral characteristics}

\maketitle

\vspace{2ex}
\textbf{\bf AMS 2010 subject classifications.} Primary: 60G10, 60G25, 60G35, Secondary: 62M20, 62P20, 93E10, 93E11

\theoremstyle{plain}

\section{Introduction}

The non-stationary and long memory time series models are of constant interest of researchers in the past decade (see, for example, papers by Dudek, Hurd and Wojtowicz [5], Johansen and Nielsen [12], Reisen et al. [32]).
These models are used when analyzing data which arise in different field of economics, finance, climatology, air pollution, signal processing.

Since the first edition of the book by Box and Jenkins (1970), autoregressive moving average (ARMA) models integrated of order $d$ are a standard tool for time series
analysis. These models are described by the equation
\be
 \psi (B) (1-B) ^ dx_t = \theta (B) \eps_t. \label{ARIMA_model} \ee
 where $ \eps_t $, $ t \in \mr Z $, are zero mean i.i.d. random variables,  $ \psi (z) $, $ \theta ( z) $ are polynomials of $ p $ and $ q $ degrees
 respectively with roots  outside the unit circle.
  This integrated ARIMA model is generalized by adding a seasonal component. A new model is described by the equation (see new edition of the book by Box and Jenkins [3] for detailes)
 \be
 \Psi (B ^ s) (1-B ^ s) ^ Dx_t = \Theta (B ^ s) \eps_t, \label{seasonal_2_model} \ee
where $ \Psi (z) $ and $ \Theta (z) $ are polynomials of degrees of $ P $ and $ Q $ respectively which have roots outside the unit circle.

When the ARIMA  sequence determined by equation (1) is inserted into (2) instead of $ \eps_t$ we have
a  general multiplicative model
\be
 \Psi (B ^ s) \psi (B) (1-B) ^ d (1-B ^ s) ^ Dx_t = \Theta (B ^ s) \theta (B) \eps_t \label {seasonal_3_model} \ee
with parameters $ (p, d, q) \times (P, D, Q)_s $, $d,D\in\mr N^*$, called SARIMA $(p, d, q)\times(P, D, Q)_s$ model.

A good performance is shown by models which include a fractional integration, that is when parameters $d$ and $D$ are fractional.
We refer to the paper by Porter-Hudak [31] who studied a seasonal ARFIMA model and applied it to the monetary aggregates used by U.S. Federal Reserve.

Another type of non-stationarity is described by periodically correlated, or cyclostationary, processes introduced by Gladyshev [8].
These processes are widely used in signal processing and communications (see Napolitano [29] for a review of the recent works on cyclostationarity and its applications).
Periodic time series may be considered as an extension of a SARIMA model (see Lund [18] for a test assessing if a PARMA model is preferable to a SARMA one) and are suitable for forecasting stream flows with quarterly, monthly or weekly cycles (see Osborn [30]).
Baek, Davis and Pipiras [1]  introduced a periodic dynamic factor model (PDFM) with periodic vector autoregressive (PVAR) factors, in contrast to seasonal VARIMA factors.
Basawa, Lund and Shao [2] investigated first-order seasonal autoregressive processes with periodically varying parameters.

The models mentioned above are used in estimation of model's parameters and forecast issues.
Note, that direct application of the developed results to real data may lead to significant increasing of values of errors of estimates due to the presence of outliers, measurement errors, incomplete information about the spectral, or model, structure etc.
This is a reason of increasing interest to robust methods of estimation that are reasonable in such cases.
For example, Reisen  et al. [33]  proposed a semiparametric robust estimator for the fractional parameters in the SARFIMA model and illustrated its application
to forecasting of sulfur dioxide $SO_2$ pollutant concentrations. Solci et al. [35] proposed robust estimates of periodic  autoregressive (PAR) model.

Robust approaches  are successfully  applied  to the problem of estimation of linear functionals from unobserved  values of stochastic processes.
The paper by Grenander [9] should be marked as the first one where the minimax extrapolation problem for stationary processes was
formulated as a game of two players and solved.
 Hosoya [11], Kassam  [14], Kassam and Poor [15], Franke [6], Vastola and  Poor [36], Moklyachuk [22, 23] studied minimax  extrapolation (forecasting), interpolation (missing values estimation) and filtering (smoothing) problems for the stationary sequences and processes. Recent results of minimax extrapolation problems for stationary vector processes and periodically correlated processes belong to  Moklyachuk and Masyutka [25, 26] and  Moklyachuk and  Golichenko (Dubovetska) [4, 24] respectively. Processes with stationary increments are investigated by Luz and Moklyachuk  [19,20]. We also mention works  by Moklyachuk and Sidei [27,28], Masyutka, Moklyachuk and Sidei [21], who derive minimax estimates of stationary processes from observations with missed values.  Moklyachuk and  Kozak [17] studied the problem of interpolation of stochastic sequences with periodically stationary increments.

In this article we present results of investigation of the estimation problem for stochastic sequences with periodically stationary increments.
In Section 2 we give definition of stochastic sequences $\xi(m)$ with periodically stationary (periodically correlated) increments.
These non-stationary stochastic sequences combine  periodic structure of covariation functions of sequences as well as the integrating one.
The section also contains a short review of the spectral theory of vector-valued stationary increment sequences.
Section 3 deals with the classical estimation problem for linear functionals in the case where spectral structure of the sequences $\xi(m)$ and $\eta(m)$ are exactly  known.
Estimates are obtained by representing the sequence $\xi(m)$ with periodically stationary increments as a vector sequence $\vec\xi (m)$ with stationary  increments and applying the Hilbert space projection technique.
In Section 4, we derive the minimax-robust estimates in the case, where spectral densities of sequences are not exactly known
while some sets of admissible spectral densities are specified.
In Subsection 4.1 we describe relations which determine the least favourable spectral densities and the minimax spectral characteristics of the optimal estimates of linear functionals  for some sets of admissible spectral densities which are generalizations of the corresponding sets of admissible spectral densities described in a survey article by
Kassam and Poor [15] for the case of stationary stochastic processes.

\section{Stochastic sequences with periodically stationary increments}\label{spectral_ theory}
In this section, we present a brief review of the spectral theory of stochastic sequences with periodically stationary $n$th increments.

 Consider a stochastic sequence $\{\xi(m),m\in\mathbb Z\}$. By $B_{\mu}$ denote a backward shift operator with the step $\mu\in
\mathbb Z$, such that $B_{\mu}\xi(m)=\xi(m-\mu)$; $B:=B_1$.
Recall the following definition [20, 37].

\begin{ozn}
For a given stochastic sequence $\{\xi(m),m\in\mathbb Z\}$, the
sequence
\begin{equation}
\label{Pryrist}
\xi^{(n)}(m,\mu)=(1-B_{\mu})^n\xi(m)=
\\
=\sum_{l=0}^n(-1)^l{n \choose l}\xi(m-l\mu),
\end{equation}
where ${n \choose l}=\frac{n!}{l!(n-l)!}$, is called
stochastic $n$th increment sequence with step $\mu\in\mathbb Z$.
\end{ozn}

The stochastic $n$th increment sequence $\xi^{(n)}(m,\mu)$ satisfies the
following relations:
\[ \xi^{(n)}(m,-\mu)=(-1)^{n}\xi^{(n)}(m+n\mu,\mu),\]
\[ \xi^{(n)}(m,k\mu)=\sum_{l=0}^{(k-1)n}A_{l}\xi^{(n)}(m-l\mu,\mu), \quad k\in\mathbb{N}, \]
  where coefficients $\{A_l,l=0,1,2,\ldots,(k-1)n\}$ are
determined by the repre\-sen\-tation
\[(1+x+\ldots+x^{k-1})^n=\sum_{l=0}^{(k-1)n
}A_lx^l.
\]

\begin{ozn}
\label{oznStPryrostu}
The stochastic $n$th increment sequence $\xi^{(n)}(m,\mu)$ generated
by  stochastic sequence $\{\xi(m),m\in\mathbb Z\}$ is wide sense
stationary if the mathematical expectations
\[E\xi^{(n)}(m_0,\mu)=c^{(n)}(\mu),\]
\[E\xi^{(n)}(m_0+m,\mu_1)\xi^{(n)}(m_0,\mu_2)=D^{(n)}(m,\mu_1,\mu_2)\]
exist for all $m_0,\mu,m,\mu_1,\mu_2$ and do not depend on $m_0$.
The function $c^{(n)}(\mu)$ is called the mean value of the $n$th
increment sequence $\xi^{(n)}(m,\mu)$ and the function $D^{(n)}(m,\mu_1,\mu_2)$ is
called the structural function of the stationary $n$th increment
sequence (or structural function of $n$th order of the stochastic
sequence $\{\xi(m),m\in\mathbb Z\}$).
\\
The stochastic sequence $\{\xi(m),m\in\mathbb   Z\}$ which
determines the stationary $n$th increment sequence
$\xi^{(n)}(m,\mu)$ by formula (4) is called the stochastic
sequence with stationary $n$th increments (or integrated sequence of order $n$)
\end{ozn}

\begin{thm}\label{thm1}
The mean value $c^{(n)}(\mu)$ and the structural function
$D^{(n)}(m,\mu_1,\mu_2)$ of the stochastic stationary $n$th
increment sequence $\xi^{(n)}(m,\mu)$ can be represented in the
 forms
\begin{equation}
\label{serFnaR}
c^{(n)}(\mu)=c\mu^n,
\end{equation}
\begin{equation}
\label{strFnaR1}
 D^{(n)}(m;\mu_1,\mu_2)=
 \\
 \int_{-\pi}^{\pi}e^{i\lambda
m} (1-e^{-i\mu_1\lambda})^n(1-e^{i\mu_2\lambda})^n\frac{1}
{\lambda^{2n}}dF(\lambda),
\end{equation}
where $c$ is a constant, $F(\lambda)$ is a left-continuous
nondecreasing bounded function such that $F(-\pi)=0$. The constant $c$
and the function $F(\lambda)$ are determined uniquely by the
increment sequence $\xi^{(n)}(m,\mu)$.
\\
On the other hand, a function $c^{(n)}(\mu)$ which has form
$(5)$ with a constant $c$ and a function
$D^{(n)}(m;\mu_1,\mu_2)$ which has form $(6)$ with a
function $F(\lambda)$ which satisfies the indicated conditions are
the mean value and the structural function of a stationary $n$th
increment sequence $\xi^{(n)}(m,\mu)$.
\end{thm}

Note that we will call by spectral function and spectral density of the stochastic sequence with stationary increments the spectral function and the spectral density of the corresponding stationary increment sequence.

Making use of representation $(6)$ and the Karhunen theorem [7,13] one can obtain the spectral representation of
the stationary $n$th increment sequence $\xi^{(n)}(m,\mu)$:
\begin{equation}
\label{predZnaR}
\xi^{(n)}(m,\mu)=\int_{-\pi}^{\pi}
e^{im\lambda}(1-e^{-i\mu\lambda})^n\frac{1}{(i\lambda)^n}dZ_{\xi^{(n)}}(\lambda),
\end{equation}
where $Z_{\eta^{(n)}}(\lambda)$ is a stochastic process with uncorrelated increments on $[-\pi,\pi)$ connected with the spectral function $F(\lambda)$ by
the relation $(-\pi\leq \lambda_1<\lambda_2<\pi)$
\[
 \mt E|Z_{\xi^{(n)}}(\lambda_2)-Z_{\xi^{(n)}}(\lambda_1)|^2=F(\lambda_2)-F(\lambda_1)<\infty.
\]

\begin{ozn}
\label{OznPeriodProc}
A stochastic sequence $\{\xi(m),m\in\mathbb Z\}$ is called  stochastic
sequence with periodically stationary (periodically correlated) increments with period $T$ if the  $n$th increment sequence
\begin{equation*}
\xi^{(n)}(m,\mu T)=(1-B_{\mu T})^n\xi(m)
\end{equation*}
is stationary.
\end{ozn}

It follows from  Definition 2.3 that the sequence
\begin{equation}
\label{PerehidXi}
\xi_{p}(m)=\xi(mT+p-1),  p=1,2,\dots,T; m\in\mathbb Z
\end{equation}
forms a vector-valued sequence
$\vec{\xi}(m)=\left\{\xi_{p}(m)\right\}_{p=1,2,\dots,T}, m\in\mathbb Z$
with stationary $n$th increments. Really, for all $p=1,2,\dots,T$,
\begin{multline*}
\xi_{p}^{(n)}(m,\mu)=\sum_{l=0}^n(-1)^l \binom{n}{l}\xi_{p}(m-l\mu)=\\
=\sum_{l=0}^n(-1)^l \binom{n}{l}\xi((m-l\mu)T+p-1)=\xi^{(n)}(mT+p-1,\mu T),
\end{multline*}
where $\xi_{p}^{(n)}(m,\mu)$ is the $n$th increment of the $p$-th component of the vector-valued sequence $\vec{\xi}(m)$.

\begin{thm}\label{thm1}
The structural function
$D^{(n)}(m,\mu_1,\mu_2)$ of the vector-valued stochastic stationary $n$th
increment sequence $\vec{\xi}^{(n)}(m,\mu)$ can be represented in the form
\begin{equation}
\label{strFnaR}
 D^{(n)}(m;\mu_1,\mu_2)=
 \\
 =\int_{-\pi}^{\pi}e^{i\lambda
m} (1-e^{-i\mu_1\lambda})^n(1-e^{i\mu_2\lambda})^n\frac{1}
{\lambda^{2n}}dF(\lambda),
\end{equation}
where $F(\lambda)$ is the matrix-valued spectral function of the stationary stochastic sequence $\vec{\xi}^{(n)}(m,\mu)$.
\\
The stationary $n$th increment sequence $\vec{\xi}^{(n)}(m,\mu)$ admits the spectral representation
\begin{equation}
\label{SpectrPred}
\vec{\xi}^{(n)}(m,\mu)=\int_{-\pi}^{\pi}e^{im\lambda}(1-e^{-i\mu\lambda})^{n}\frac{1}{(i\lambda)^{n}}d\vec{Z}_{\xi^{(n)}}(\lambda),
\end{equation}
where $d\vec{Z}_{\xi^{(n)}}(\lambda)=\{Z_{ p}(\lambda)\}_{p=1}^{T}$ is a (vector-valued) stochastic process with uncorrelated increments on $[-\pi,\pi)$ connected with the spectral function $F(\lambda)$ by
the relation
\begin{equation*}
 \mt E(Z_{p}(\lambda_2)-Z_{p}(\lambda_1))(\overline{ Z_{q}(\lambda_2)-Z_{q}(\lambda_1)})= F_{pq}(\lambda_2)-F_{pq}(\lambda_1),  -\pi\leq \lambda_1<\lambda_2<\pi.
\end{equation*}
\end{thm}

\section{Hilbert space projection method of estimation}\label{classical_estimstion}

Consider a vector-valued stochastic sequence with stationary $n$th increments $\vec{\xi}(m)$ constructed from the sequence $\xi(m)$ with the help of transformation (8).
Let the stationary $n$th increment sequence $\vec{\xi}^{(n)}(m,\mu)=\{\xi_{p}^{(n)}(m,\mu)\}_{p=1}^{T}$
has an absolutely continuous spectral function $F(\lambda)$
and the spectral density $f(\lambda)=\{f_{ij}(\lambda)\}_{i,j=1}^{T}$.

Let  $\vec{\eta}(m)=\{\eta_{p}(m)\}_{p=1}^{T}$ be an uncorrelated with
the sequence $\xi(m)$ stationary stochastic sequence with absolutely
continuous spectral function $G(\lambda)$ and spectral density
$g(\lambda)=\{g_{ij}(\lambda)\}_{i,j=1}^{T}$.

We will assume that the mean
values of the increment sequence $\vec{\xi}^{(n)}(m,\mu)$ and stationary
sequence  $\vec{\eta}(m)$ equal to 0. We will also consider the increment step $\mu>0$.

Consider the problem of mean square optimal linear estimation of the functional
\begin{equation}
A_{N}\vec{\xi}=\sum_{k=0}^{N}(\vec{a}(k))^{\top}\vec{\xi}(k)
\end{equation}
which depend on the unobserved values of the stochastic sequence $\vec{\xi}(k)=\{\xi_{p}(k)\}_{p=1}^{T}$ with stationary $n$th
increments.
Estimates are based on observations of the sequence $\vec\zeta(m)=\vec\xi(m)+\vec\eta(m)$ at points of the set
$\mr Z\setminus\{0,1,2\ldots,N\}$.

Assume that spectral densities $f(\lambda)$ and $g(\lambda)$ satisfy the minimality condition
\be
 \ip \text{Tr}\left[ \frac{\lambda^{2n}}{|1-e^{i\lambda\mu}|^{2n}}(f(\lambda)+{\lambda}^{2n}g(\lambda))^{-1}\right]
 d\lambda<\infty.
\label{umova11_e_d}
\ee
This is the necessary and sufficient condition under which the mean square errors of the optimal estimates of the functional $A_N\vec\xi$ is not equal to $0$.

The classical Hilbert space estimation technique proposed by Kolmogorov [16] can be described as a $3$-stage procedure:
(i) define a target element of the space $H=L_2(\Omega, \mathcal{F},\mt P)$ to be estimated,
(ii) define a subspace of $H$ generated by observations,
(iii) find an estimate of the target element as an orthogonal projection on the defined subspace.

\textbf{Stage i}. The  functional $A_{N}\vec{\xi}$ does not  belong to the space $H$.
With the help of the following lemma we  describe representations of the functional   as a sum of a functional with finite second moments  belonging to $H$ and a functional depending on observed values of the sequence $\vec\zeta(k)$ (``initial values'').

\begin{lema}\label{lema predst A}
The functional $A_N\vec\zeta$ admits the representation
\be \label{zobrazh A_N_i_st.n_d}
    A_N\vec\xi=A_N\vec\zeta-A_N\vec\eta=H_N\vec\xi-V_N\vec\zeta,
\ee
where
\[
    H_N\vec\xi:=B_N\vec\zeta-A_N\vec\eta,\]
    \[
    A_{N}\vec{\zeta}=\sum_{k=0}^{N}(\vec{a}(k))^{\top}\vec{\zeta}(k),\quad
    A_{N}\vec{\eta}=\sum_{k=0}^{N}(\vec{a}(k))^{\top}\vec{\eta}(k),\]
\[
B_{N}\vec{\zeta}=\sum_{k=0}^{N}(\vec{b}_{N}(k))^{\top}\vec{\zeta}^{(n)}(k,\mu),\quad
V_{N}\vec{\zeta}=\sum_{k=-\mu n}^{-1} (\vec{v}_{N}(k))^{\top}\vec{\zeta}(k),
\]
the coefficients
 $\vec{v}_{N}(k))=\{ v_{N,p}(k)\}_{p=1}^{T}, k=-1,-2,\dots,-\mu n$ and $\vec{b}_{N}(k))=\{b_{N,p}(k)\}_{p=1}^{T}, k=0,1,\dots,N$
are calculated by the formulas
\begin{equation}
\label{koefv_N_diskr}
 v_{N,p}(k)=\sum_{l=\left[-\frac{k}{\mu}\right]'}^{\min\left\{\left[\frac{N-k}{\mu}\right],n\right\}}(-1)^{l} \binom{n}{l} b_{N,p}(l\mu+k), \\
  k=-1,-2,\dots,-\mu n,
\end{equation}
\begin{equation} \label{determ_b_N}
b_{N,p}(k)=\sum_{m=k}^{N}a_{p}(m)d_{\mu}(m-k)=(D^{\mu}_{N}{\me a}_{N,p})_k, \\ k=0,1,\dots,N,
\end{equation}
coefficients $\{d_{\mu}(k):k\geq0\}$ are determined by the relationship
\[
 \sum_{k=0}^{\infty}d_{\mu}(k)x^k=\left(\sum_{j=0}^{\infty}x^{\mu j}\right)^n,\]
$D^{\mu}_{N}$ is a  linear transformation  determined by a matrix with the entries $(D^{\mu}_{N})_{k,j}=d_{\mu}(j-k)$ if
$0\leq k\leq j\leq N$, and $(D^{\mu}_{N})_{k,j}=0$ if $0\leq j<k\leq N$;
$D^{\mu}_N{\me a}_N=\{D^{\mu}_{N}{\me a}_{N,p}\}_{p=1}^T$,
${\me a}_{N,p}=(a_p(0),a_p(1),a_p(2),\ldots,a_p(N))^{\top}$, $p=1,2,\dots,T$.
\end{lema}

The functional $H_N\vec\xi$ from representation (13) has finite variance and the functional  $V_N\vec\zeta$ depends on the known observations of the stochastic sequence $\vec\zeta(k)$ at points $k=-\mu n,-\mu n+1,\ldots,-1$. Therefore, estimates $\widehat{A}_N\vec\xi$ and $\widehat{H}_N\vec\xi$ of the functionals $A_N\vec\xi$ and $H_N\vec\xi$ and the mean-square errors $\Delta(f,g;\widehat{A}_N\vec\xi)=\mt E |A_N\vec\xi-\widehat{A}_N\vec\xi|^2$ and $\Delta(f,g;\widehat{H}_N\vec\xi)=\mt E
|H_N\vec\xi-\widehat{H}_N\vec\xi|^2$ of the estimates $\widehat{A}_N\vec\xi$ and $\widehat{H}_N\vec\xi$ satisfy the following relations
\be\label{mainformula}
    \widehat{A}_N\vec\xi=\widehat{H}_N\vec\xi-V_N\vec\zeta,\ee
\begin{equation*}
    \Delta(f,g;\widehat{A}_N\vec\xi)
    =\mt E |A_N\vec\xi-\widehat{A}_N\vec\xi|^2=
      \\
      \mt E|H_N\vec\xi-\widehat{H}_N\vec\xi|^2=\Delta(f,g;\widehat{H}_N\vec\xi).
    \end{equation*}
Therefore, the estimation problem for the functional $A_N\vec\xi$ is equivalent to the estimation problem for the functional $H_N\vec\xi$.
This problem can be solved by applying the Hilbert space projection method proposed by Kolmogorov [16].

The stationary stochastic sequence $\vec\eta(m)$ admits the spectral representation
\[\vec\eta(m)=\int_{-\pi}^{\pi}e^{i\lambda m}d\vec{Z}_{\eta}(\lambda),\]
where $\vec{Z}_{\eta}(\lambda)$   is a random process with uncorrelated increments on
$[-\pi,\pi)$ corresponding to the spectral function $G(\lambda)$.
The random processes $\vec{Z}_{\eta}(\lambda)$ and $\vec{Z}_{\eta^{(n)}}(\lambda)$ are connected by the relation $d\vec{Z}_{\eta^{(n)}}(\lambda)=(i\lambda)^n d\vec{Z}_{\eta }(\lambda)$,
$\lambda\in[-\pi,\pi)$, obtained in [19]. The spectral density $p(\lambda)$ of the
sequence $\vec\zeta(m)$ is determined by spectral densities $f(\lambda)$
and $g(\lambda)$ by the relation
\[p(\lambda)=f(\lambda)+\lambda^{2n}g(\lambda).\]

With the help of the spectral representations of stochastic sequences involved we can write the following spectral representation of the  functional
\begin{equation*}
H_N\vec\xi=\\
\int_{-\pi}^{\pi}\left(\vec{B}_{\mu,N}(e^{i\lambda})\right)^{\top}\frac{(1-e^{-i\lambda \mu})^n}{(i\lambda)^n}d\vec{Z}_{\xi^{(n)}+\eta^{(n)}}(\lambda)-
\\
-\int_{-\pi}^{\pi}\left(\vec{A}_N(e^{i\lambda})\right)^{\top}d\vec{Z}_{\eta}(\lambda),
\end{equation*}
where
\begin{equation*}
\vec{B}_{\mu,N}(e^{i\lambda})=\sum_{k=0}^{N}\vec{b}_{\mu,N}(k)e^{i\lambda k}=\\
=\sum_{k=0}^{N}(D^{\mu}_N\me a_N)_ke^{i\lambda k},
 \vec{A}_N(e^{i\lambda})=\sum_{k=0}^{N}\vec{a}(k)e^{i\lambda k}.
\end{equation*}

At \textbf{stage ii}, we deal with the following notations. Denote by $H^{0-}(\xi^{(n)}_{\mu}+\eta^{(n)}_{\mu})$ the closed linear subspace generated by
values $\{\xi_p^{(n)}(k,\mu)+\eta_p^{(n)}(k,\mu):p=1,\dots,T;k=-1,-2,-3,\dots\}$
in the Hilbert space $H=L_2(\Omega,\mathcal{F},\mt P)$ of random variables $\gamma$ with zero mean value, $\mt E\gamma=0$, finite variance, $\mt E|\gamma|^2<\infty$, and the inner product $(\gamma_1;\gamma_2)=\mt E\gamma_1\overline{\gamma_2}$.

Denote by $H^{N+}(\xi^{(n)}_{-\mu}+\eta^{(n)}_{-\mu})$ the closed linear subspace of the Hilbert space $H=L_2(\Omega,\mathcal{F},\mt P)$  generated by elements
$\{\xi^{(n)}_p(k,-\mu)+\eta^{(n)}_p(k,-\mu):p=1,\dots,T;k\geq N+1\}$.

The equality $\xi^{(n)}_p(k,-\mu)=(-1)^n\xi^{(n)}_p(k+\mu n,\mu)$ implies
\[
    H^{N+}(\xi_{-\mu}^{(n)}+\eta_{-\mu}^{(n)})=H^{(N+\mu n)+}(\xi_{\mu}^{(n)}+\eta_{\mu}^{(n)}).\]

Denote by $L_2^{0-}(f(\lambda)+\lambda^{2n}g(\lambda))$ and $L_2^{N+}(f(\lambda)+\lambda^{2n}g(\lambda))$ the closed linear subspaces of the Hilbert space
$L_2(f(\lambda)+\lambda^{2n}g(\lambda))$  of vector-valued functions with the inner product $\langle g_1;g_2\rangle=\ip (g_1(\lambda))^{\top}(f(\lambda)+\lambda^{2n}g(\lambda))\overline{g_2(\lambda)}d\lambda$ which is generated by the functions
\begin{equation*}
e^{i\lambda k}(1-e^{-i\lambda \mu})^n\frac{1}{(i\lambda)^n}{\delta}_{l},\quad  {\delta}_l=\{\delta_{lp}\}_{p=1}^{T}, \\
 l=1,\dots,T;\,\, k \leq -1,
;\,\, \,k\geq N+1,
\end{equation*}
respectively, where $\delta_{lp}$ are Kronecker symbols.

The representation
\begin{equation*}
 \vec\xi^{(n)}(k,\mu)+\vec\eta^{(n)}(k,\mu)=\int_{-\pi}^{\pi} e^{i\lambda
k}(1-e^{-i\lambda\mu})^n\dfrac{1}{(i\lambda)^{n}}d\vec{Z}_{\xi^{(n)}+\eta^{(n)}}(\lambda)
\end{equation*}
yields a  one to one correspondence between elements $e^{i\lambda k}(1-e^{-i\lambda\mu})^n(i\lambda)^{-n}$ of the space
$$L_2^{0-}(f(\lambda)+\lambda^{2n}g(\lambda))\oplus L_2^{N+}(f(\lambda)+\lambda^{2n}g(\lambda))$$
and elements
$\vec\xi^{(n)}(k,\mu)+\vec\eta^{(n)}(k,\mu)$
of the space
\begin{equation*}
H^{0-}(\xi^{(n)}_{\mu}+\eta^{(n)}_{\mu})\oplus H^{N+}(\xi^{(n)}_{-\mu}+\eta^{(n)}_{-\mu})=
H^{0-}(\xi^{(n)}_{\mu}+\eta^{(n)}_{\mu})\oplus H^{(N+\mu n)+}(\xi^{(n)}_{\mu}+\eta^{(n)}_{\mu}).
\end{equation*}

Relation (16) implies that every linear estimate $\widehat{A}\vec\xi$ of the functional $A\vec\xi$
can be represented in the form
\begin{equation}
 \label{otsinka A_e_d}
 \widehat{A}_N\vec\xi=\ip
(\vec{h}_{\mu,N}(\lambda))^{\top}d\vec{Z}_{\xi^{(n)}+\eta^{(n)}}(\lambda)-\sum_{k=-\mu n}^{-1}(\vec v_{\mu}(k))^{\top}(\vec\xi(k)+\vec\eta(k)),
\end{equation}
 where
$\vec{h}_{\mu,N}(\lambda)=\{h_{p}(\lambda)\}_{p=1}^{T}$ is the spectral characteristic of the optimal estimate $\widehat{H}_N\vec\xi$.

At \textbf{stage iii}, we find the mean square optimal estimate
$\widehat{H}_N\vec\xi$ as a projection of the element $H_N\vec\xi$ on the
subspace $H^{0-}(\xi^{(n)}_{\mu}+\eta^{(n)}_{\mu})\oplus H^{(N+\mu n)+}(\xi^{(n)}_{\mu}+\eta^{(n)}_{\mu})$. This projection is
determined by two conditions:

1) $ \widehat{H}_N\vec\xi\in H^{0-}(\xi^{(n)}_{\mu}+\eta^{(n)}_{\mu})\oplus H^{(N+\mu n)+}(\xi^{(n)}_{\mu}+\eta^{(n)}_{\mu}) $;

2) $(H_N\vec\xi-\widehat{H}_N\vec\xi)
\perp
H^{0-}(\xi^{(n)}_{\mu}+\eta^{(n)}_{\mu})\oplus H^{(N+\mu n)+}(\xi^{(n)}_{\mu}+\eta^{(n)}_{\mu})$.

The second condition implies the following relation which holds true for all $k\leq-1$ and $k\geq N+\mu n+1$
\bml
\int_{-\pi}^{\pi}
\bigg[\bigg(\vec{B}_{\mu,N}(e^{i\lambda})\frac{(1-e^{-i\lambda\mu})^n}{(i\lambda)^n}-
\vec{h}_{\mu,N}(\lambda)\bigg)^{\top}
p(\lambda)-
\\
-(\vec{A}_N(e^{i\lambda})^{\top}g(\lambda)(-i\lambda)^n\bigg]
\frac{(1-e^{i\lambda\mu})^n}{(-i\lambda)^{n}}e^{-i\lambda k}d\lambda=0.
 \end{multline*}

\noindent This relation allows us to derive the spectral characteristic
$\vec{h}_{\mu,N}(\lambda)$ of the estimate $\widehat{H}_N\vec\xi$ which can be represented in the form
\begin{multline} \label{spectr A}
(\vec{h}_{\mu,N}(\lambda))^{\top}=(\vec{B}_{\mu,N}(e^{i\lambda}))^{\top}
\frac{(1-e^{-i\lambda \mu})^n}{(i\lambda)^n}-
\\
-
\vec{A}_N(e^{i\lambda})^{\top}g(\lambda)(-i\lambda)^n
(p(\lambda))^{-1}
-
\frac{(-i\lambda)^{n}(\vec{C}_{\mu,N}(e^{i\lambda}))^{\top}}{(1-e^{i\lambda \mu})^n}(p(\lambda))^{-1},
\end{multline}
\[
\vec{C}_{\mu,N}(e^{i \lambda})=\sum_{k=0}^{N+\mu N}\vec{c}_{\mu,N}(k)e^{ik\lambda},\]
 $\vec{c}_{\mu,N}(k)=\{c_{\mu,N,p}(k)\}_{p=1}^T, k=0,1,\dots,N+\mu n,$  are unknown coefficients to be found.

It follows from condition 1) that the following equations should be satisfied  for $0\leq j\leq N+\mu n$
\begin{multline}
\label{eq_C}
\int_{-\pi}^{\pi} \biggl[(\vec{B}_{\mu,N}(e^{i\lambda}))^{\top}-
(\vec{A}_N(e^{i\lambda})^{\top}g(\lambda)
\frac{(\lambda)^{2n}}{(1-e^{-i\lambda \mu})^n}
(p(\lambda))^{-1}-
\\
-
\frac{\lambda^{2n}(\vec{C}_{\mu,N}(e^{i\lambda}))^{\top}}{(1-e^{-i\lambda \mu})^n (1-e^{i\lambda \mu})^n}(p(\lambda))^{-1}\biggr]e^{-ij\lambda}d\lambda=0.
\end{multline}

 Define for $0\leq k,j \leq N$ the Fourier coefficients of the corresponding functions
\[
T^{\mu}_{k,j}=\frac{1}{2\pi}\int_{-\pi}^{\pi}
e^{i\lambda(j-k)}\frac{\lambda^{2n}g(\lambda)}{|1-e^{i\lambda\mu}|^{2n}}
(p(\lambda))^{-1}
d\lambda;
\]
\[
P_{k,j}^{\mu}=\frac{1}{2\pi}\int_{-\pi}^{\pi} e^{i\lambda (j-k)}
\dfrac{\lambda^{2n}}{|1-e^{i\lambda\mu}|^{2n}}(p(\lambda))^{-1}
d\lambda;
\]
 \[
 Q_{k,j}=\frac{1}{2\pi}\int_{-\pi}^{\pi}
e^{i\lambda(j-k)}f(\lambda)g(\lambda)(p(\lambda))^{-1}d\lambda.
\]

  \noindent Making use of the defined Fourier coefficients, relation (19) can be presented as a system of $N+\mu n+1$ linear equations
  determining the unknown coefficients ${\vec c}_{\mu,N}(k)$, $0\leq k\leq N+\mu n$.
\begin{equation}
\label{linear equations1}
    \vec{b}_{\mu,N}(j)-\sum_{ m=0}^{N+\mu n}T^{\mu}_{j,m}\vec{a}_{\mu,N}(m)
    =\sum_{k=0}^{N+\mu n}P_{j,k}^{\mu}\vec{c}_{\mu,N}(k), 0\leq j\leq N,
    \end{equation}
\begin{equation}
\label{linear equations2}
    -\sum_{ m=0}^{N+\mu n}T^{\mu}_{j,m}\vec{a}_{\mu,N}(m)
    =\sum_{k=0}^{N+\mu n}P_{j,k}^{\mu}\vec{c}_{\mu,N }(k), N+1\leq j\leq N+\mu n,
    \end{equation}
where coefficients $\{\vec{a}_{\mu,N}(m):0\leq m\leq N+\mu n\}$ are calculated by the formula
\begin{equation}
\label{coeff a_N_mu}
    \vec{a}_{\mu,N}(m)
    =\sum_{l=\max\ld\{\ld[\frac{m-N}{\mu}\rd]',0\rd\}}^{\min\ld\{\ld[\frac{m}{\mu}\rd],n\rd\}}
    (-1)^l{n \choose l}\vec{a}(m-\mu l), 0\leq m\leq N+\mu n,
    \end{equation}

  \noindent Denote by $[D_N^{\mu}\me a_N]_{+\mu n}$ a vector of dimension  $(N+\mu n +1)T$ which is constructed by adding
$(\mu n)T$ zeros to the vector $D_N^{\mu}\me a_N$ of dimension  $(N+1)T$.
Making use of this definition the system (20) -- (21) can be represented in the matrix form:
\[[D_N^{\mu}\me a_N]_{+\mu n}-\me T^{\mu}_N\me a^{\mu}_N=\me P^{\mu}_N\me c^{\mu}_N,\]
\begin{equation*}
 \me a^{\mu}_N=((\vec{a}_{\mu,N}(0))^{\top},(\vec{a}_{\mu,N}(1))^{\top},
 \ldots,(\vec{a}_{\mu,N}(N+\mu n))^{\top})^{\top}
\end{equation*}
\begin{equation*}
    \me c^{\mu}_N=((\vec{c}_{\mu,N}(0))^{\top},(\vec{c}_{\mu,N}(1))^{\top},
    \ldots,(\vec{c}_{\mu,N}(N+\mu n))^{\top})^{\top}
\end{equation*}
are  vectors of  dimension $(N+\mu n+1)T$; $\me P^{\mu}_N$ and $\me T^{\mu}_N$ are matrices of dimension $(N+\mu n+1)T\times(N+\mu n+1)T$ with $T\times T$ matrix elements $(\me P^{\mu}_N)_{j,k}=P_{j,k}^{\mu}$ and $(\me T^{\mu}_N)_{j, k} =T^{\mu}_{j,k}$, $0\leq j,k\leq N+\mu n$.

Thus, the coefficients $\vec{c}_{\mu,N}(k)$, $0\leq k\leq N+\mu n$, are determined by the formula $ 0\leq k\leq N+\mu n$
\[
  \vec{c}_{\mu,N}(k)=\ld((\me P^{\mu}_N)^{-1}[D_N^{\mu}\me a_N]_{+\mu n}-(\me P^{\mu}_N)^{-1}\me T^{\mu}_N\me a^{\mu}_{N}\rd)_k,\]
where $\ld((\me P^{\mu}_N)^{-1}[D_N^{\mu}\me a_N]_{+\mu n}-(\me P^{\mu}_N)^{-1}\me T^{\mu}_N\me a^{\mu}_N\rd)_k$, $0\leq k\leq N+\mu n$, is the
 $k$th element of the vector  $(\me P^{\mu}_N)^{-1}[D_N^{\mu}\me a_N]_{+\mu n}-(\me P^{\mu}_N)^{-1}\me T^{\mu}_N\me a^{\mu}_N$.

The existence of the inverse matrix $(\me P^{\mu}_N)^{-1}$ was shown in [20] under condition (12).

The spectral characteristic  $\vec{h}_{\mu,N}(\lambda)$ of the estimate $\widehat{H}_N\xi$ of the functional $H_N\xi$ is calculated by formula (18), where
\begin{equation}
\label{spectr C}
 \vec{C}_{\mu,N}(e^{i \lambda})=\sum_{k=0}^{N+\mu n}
    \bigg((\me P^{\mu}_N)^{-1}[D_N^{\mu}\me a_N]_{+\mu n}- (\me P^{\mu}_N)^{-1}\me T^{\mu}_N\me a^{\mu}_N\bigg)_k e^{i\lambda k}.
\end{equation}

The value of the  mean-square errors of the estimates $\widehat{A}_N\vec\xi$ and $\widehat{H}_N\vec\xi$ can be calculated by the formula

\[
\Delta(f,g;\widehat{A}\vec\xi)=\Delta(f,g;\widehat{H}\vec\xi)= \mt E|H\vec\xi-\widehat{H}\vec\xi|^2=
\]
\[
=
\frac{1}{2\pi}\int_{-\pi}^{\pi}
\frac{\lambda^{2n}}{|1-e^{i\lambda\mu}|^{2n}}
\left[{(1-e^{i\lambda\mu})^n}(\vec{A}_N(e^{i\lambda}))^{\top}g(\lambda) +
(\vec{C}_{\mu,N}(e^{i \lambda}))^{\top}
\right](f(\lambda)+{\lambda}^{2n}g(\lambda))^{-1}\,
\times
\]
\[
\times
 f(\lambda)\, (f(\lambda)+{\lambda}^{2n}g(\lambda))^{-1}
\left[{(1-e^{-i\lambda\mu})^n}{\vec{A}_N(e^{-i\lambda})}g(\lambda) +
(\vec{C}_{\mu,N}(e^{-i \lambda}))
\right]
d\lambda+
\]
\[
+\frac{1}{2\pi}\int_{-\pi}^{\pi}
\frac{1}{|1-e^{i\lambda\mu}|^{2n}}
\left[{(1-e^{i\lambda\mu})^n}(\vec{A}_N(e^{i\lambda}))^{\top}f(\lambda) -(\lambda)^{2n}
(\vec{C}_{\mu,N}(e^{i \lambda}))^{\top}
\right](f(\lambda)+{\lambda}^{2n}g(\lambda))^{-1}
\times
\]
\[
\times
g(\lambda)\,(f(\lambda)+{\lambda}^{2n}g(\lambda))^{-1}
\left[{(1-e^{-i\lambda\mu})^n}{\vec{A}_N(e^{-i\lambda})}f(\lambda) -(\lambda)^{2n}
(\vec{C}_{\mu,N}(e^{-i \lambda}))
\right]
d\lambda=
\]
 \be
= \ld\langle [D_N^{\mu}\me a_N]_{+\mu n}- \me T^{\mu}_N\me
    a^{\mu}_N,(\me P^{\mu}_N)^{-1}[D_N^{\mu}\me a_N]_{+\mu n}
    -(\me P^{\mu}_N)^{-1}\me T^{\mu}_N\me a^{\mu}_N\rd\rangle
    \\+\langle\me Q_N\me a_N,\me
    a_N\rangle,
    \label{pohybkaAN}
    \ee

\noindent where $\me Q_N$ is a matrix of the dimension $(N+1)T\times(N+1)T$ with the $T\times T$ matrix elements $(\me Q_N)_{j,k}=Q_{j,k}$, $0\leq j,k\leq N$

\begin{thm}\label{thm_intA}
Let $\{\vec\xi(m),m\in\mr Z\}$ be a stochastic sequence which defines
the stationary $n$th increment sequence $\vec\xi^{(n)}(m,\mu)$ with the
absolutely continuous spectral function $F(\lambda)$ which has
spectral density $f(\lambda)$. Let $\{\vec\eta(m),m\in\mr Z\}$ be an
uncorrelated with the sequence $\vec\xi(m)$ stationary stochastic
sequence with an absolutely continuous spectral function
$G(\lambda)$ which has spectral density $g(\lambda)$. Let the minimality condition
(12) be satisfied. The
optimal linear estimate $\widehat{A}_N\vec\xi$ of the functional
$A_N\vec\xi$ which depends on the unknown values of elements $\vec\xi(k)$, $k=0,1,2,\ldots,N$,   from
observations of the sequence $\vec\xi (m )+\vec\eta (m )$ at points of the set $Z\setminus\{0,1,2,\ldots,N\}$
is calculated by  formula (17).
The spectral characteristic $\vec{h}_{\mu,N}(\lambda)$ of the optimal
estimate $\widehat{A}_N\vec\xi$ is calculated by  formulas (18), (23). The value of the mean-square error
$\Delta(f,g;\widehat{A}_N\vec\xi)$ is calculated by  formula (24).
\end{thm}

\begin{nas}
The spectral characteristic $\vec{h}_{\mu,N}(\lambda)$ (18) admits the representation
$\vec{h}_{\mu,N}(\lambda)=\vec{h}_{\mu,N}^{1}(\lambda)-\vec{h}_{\mu,N}^2(\lambda)$, where
\begin{multline}
\label{spectr h1}
(\vec{h}_{\mu}^1(\lambda))^{\top}=
(\vec{B}_{\mu,N}(e^{i\lambda}))^{\top}
\frac{(1-e^{-i\lambda \mu})^n}{(i\lambda)^n}-
\frac{(-i\lambda)^{n}}{(1-e^{i\lambda \mu})^n}
\times
\\
\times
\left(
\sum_{k=0}^{N+\mu n}
    \ld((\me P^{\mu}_N)^{-1}[D_N^{\mu}\me a_N]_{+\mu n}\rd)_k e^{i\lambda k}
\right)^{\top}
(p(\lambda))^{-1},
 \end{multline}
 \begin{multline} \label{spectr h2}
(\vec{h}_{\mu,N}^2(\lambda))^{\top}=
(\vec{A}_N(e^{i\lambda }))^{\top}
 {(-i\lambda)^ng(\lambda)}(p(\lambda))^{-1}-\frac{(-i\lambda)^{n}}{(1-e^{i\lambda \mu})^n}
 \times
\\
\times
\left(
\sum_{k=0}^{N+\mu n}
    \ld((\me P^{\mu}_N)^{-1}\me T^{\mu}_N\me a^{\mu}_N\rd)_k e^{i\lambda k}
\right)^{\top}
(p(\lambda))^{-1}.
\end{multline}
Here $\vec{h}_{\mu,N}^1(\lambda)$ and $\vec{h}_{\mu,N}^2(\lambda)$ are spectral characteristics of the optimal estimates
$\widehat{B}_N\vec\zeta$ and $\widehat{A}_N\vec\eta$ of the functionals $B_N\vec\zeta$ and $A_N\vec\eta$ respectively based on observations
$\vec\xi(k)+\vec\eta(k)$ at points of the set $Z\setminus\{0,1,2,\ldots,N\}$.
\end{nas}

\subsection{Estimation of stochastic sequences with periodically stationary increment}

Consider the problem of mean square optimal linear estimation of the functional
\begin{equation}
A_{M}{\vartheta}=\sum_{k=0}^{N}{a}^{(\vartheta)}(k)\vartheta(k)
\end{equation}
which depend on unobserved values of the stochastic sequence ${\vartheta}(m)$ with periodically stationary
increments. Estimates are based on observations of the sequence $\zeta(m)=\vartheta(m)+\eta(m)$ at points of the set $Z\setminus\{0,1,2,\ldots,N\}$.

The functional $A_M{\vartheta}$ can be represented in the form
\begin{multline*}
A_M{\vartheta} = \sum_{k=0}^{M}{a}^{(\vartheta)}(k)\vartheta(k)=
\\
=\sum_{m=0}^{N}\sum_{p=1}^{T}
{a}^{(\vartheta)}(mT+p-1)\vartheta(mT+p-1)=
\\
 =\sum_{m=0}^{N}\sum_{p=1}^{T}a_p(m)\xi_p(m)=\sum_{m=0}^{N}(\vec{a}(m))^{\top}\vec{\xi}(m)=A_N\vec{\xi},
\end{multline*}
where $N=[\frac{M}{T}]$, the sequence $\vec{\xi}(m) $ is determined by the formula
 \begin{equation}
 \label{zeta}
\vec{\xi}(m)=({\xi}_1(m),{\xi}_2(m),\dots,{\xi}_T(m))^{\top},\\
 {\xi}_p(m)=\vartheta(mT+p-1);\,p=1,2,\dots,T;
 \end{equation}
  \begin{multline}\label{aNzeta}
(\vec{a}(m))^{\top} = ({a}_1(m),{a}_2(m),\dots,{a}_T(m))^{\top},
  \\
  {a}_p(m) = a^{\vartheta}(mT+p-1);0\leq m\leq N; mT+p-1\leq M;
\\  {a}_p(N) = 0;
M+1\leq NT+p-1\leq (N+1)T-1.
\end{multline}

Making use of the introduced notations and statements of Theorem 3.1 we can claim that the following theorem holds true.

\begin{thm}
\label{thm_est_A_Nzeta}
Let a stochastic sequence ${\vartheta}(k)$ with periodically stationary increments generate by formula (28)  a vector-valued stochastic sequence $\vec{\xi}(m) $ which determine a
stationary stochastic $n$th increment sequence
$\vec{\xi}^{(n)}(m,\mu)$ with the spectral density matrix $f(\lambda)=\{f_{ij}(\lambda)\}_{i,j=1}^{T}$.
Let $\{\vec\eta(m),m\in\mr Z\}$, $\vec{\eta}(m)=({\eta}_1(m),{\eta}_2(m),\dots,{\eta}_T(m))^{\top},\,
 {\eta}_p(m)=\eta(mT+p-1);\,p=1,2,\dots,T, $
be uncorrelated with the sequence $\vec\xi(m)$ stationary stochastic
sequence with an absolutely continuous spectral function
$G(\lambda)$ which has spectral density $g(\lambda)$. Let the minimality condition
(12) be satisfied.
Let coefficients $\vec {a}(k), k\geqslant 0$ be determined by formula (29).
The optimal linear estimate $\widehat{A}_M\zeta$ of the functional $A_M\zeta=A_N\vec{\xi}$ based on observations of the sequence
$\zeta(m)=\vartheta(m)+\eta(m)$ at points of the set $Z\setminus\{0,1,2,\ldots,N\}$ is calculated by
 formula (17).
The spectral characteristic $\vec{h}_{\mu,N}(\lambda)=\{h_{\mu,N,p}(\lambda)\}_{p=1}^{T}$
and the value of the mean square error $\Delta(f;\widehat{A}_M\zeta)$
are calculated by formulas  (18), (23) and (24) respectively.
\end{thm}

\section{Minimax-robust method of estimation}\label{minimax_estimation}

The values of the mean square errors and the spectral characteristics of the optimal estimate
of the functional ${A}_N\vec\xi$
depending on the unobserved values of a stochastic sequence $\vec{\xi}(m)$ which determine a stationary stochastic $n$th increment sequence
$\vec{\xi}^{(n)}(m,\mu)$ with the spectral density matrix $f(\lambda)$
based on observations of the sequence
$\vec\xi(m)+\vec\eta(m)$ at points $Z\setminus\{0,1,2,\ldots,N\}$ can be calculated by formulas
(18), (23) and (24) respectively, under the condition that
spectral densities
$f(\lambda)$ and $g(\lambda)$ of stochastic sequences $\vec\xi(m)$ and
$\vec\eta(m)$ are exactly known.

In practical cases, however, spectral densities of sequences usually are not exactly known.
If in such cases a set $\md D=\md D_f\times\md D_g$ of admissible spectral densities is defined,
the minimax-robust approach to
estimation of linear functionals depending on unobserved values of stochastic sequences with stationary increments may be applied.
This method consists in finding an estimate that minimizes
the maximal values of the mean square errors for all spectral densities
from a given class $\md D=\md D_f\times\md D_g$ of admissible spectral densities
simultaneously.

To formalize this approach we present the following definitions.

\begin{ozn}
For a given class of spectral densities $\mathcal{D}=\md
D_f\times\md D_g$ the spectral densities
$f_0(\lambda)\in\mathcal{D}_f$, $g_0(\lambda)\in\md D_g$ are called
least favorable in the class $\mathcal{D}$ for the optimal linear
estimation of the functional $A_N\vec \xi$  if the following relation holds
true:
\[\Delta(f_0,g_0)=\Delta(h(f_0,g_0);f_0,g_0)=
\max_{(f,g)\in\mathcal{D}_f\times\md
D_g}\Delta(h(f,g);f,g).\]
\end{ozn}

\begin{ozn}
For a given class of spectral densities $\mathcal{D}=\md
D_f\times\md D_g$ the spectral characteristic $h^0(\lambda)$ of
the optimal linear estimate of the functional $A_N\vec \xi$ is called
minimax-robust if there are satisfied the conditions
\[h^0(\lambda)\in H_{\mathcal{D}}=\bigcap_{(f,g)\in\mathcal{D}_f\times\md D_g}L_2^{0-}(f(\lambda)+\lambda^{2n}g(\lambda)),\]
\[\min_{h\in H_{\mathcal{D}}}\max_{(f,g)\in \mathcal{D}_f\times\md D_g}\Delta(h;f,g)=\max_{(f,g)\in\mathcal{D}_f\times\md
D_g}\Delta(h^0;f,g).\]
\end{ozn}

Taking into account the introduced definitions and the derived relations we can verify that the following lemma holds true.

\begin{lema}
The spectral densities $f^0\in\mathcal{D}_f$,
$g^0\in\mathcal{D}_g$ which satisfy the minimality condition (12)
are least favorable in the class $\md D=\md D_f\times\md D_g$ for
the optimal linear estimation of the functional $A_N\vec\xi$ based on observations of the sequence $\xi(m)+\eta(m)$
at points   $m\in\mr Z\setminus\{0,1,2,\ldots,N\}$
if the matrices $ (\me P^{\mu}_N)^0$, $(\me T^{\mu}_N)^0$, $(\me Q_N)^0$ whose elements are defined by the Fourier coefficients of the functions
\[
\frac{\lambda^{2n}g^0(\lambda)}{|1-e^{i\lambda\mu}|^{2n}}
(f^0(\lambda)+{\lambda}^{2n}g^0(\lambda))^{-1},\quad
\dfrac{\lambda^{2n}}{|1-e^{i\lambda\mu}|^{2n}}(f^0(\lambda)+{\lambda}^{2n}g^0(\lambda))^{-1},
\]
\[
f^0(\lambda)g^0(\lambda)(f^0(\lambda)+{\lambda}^{2n}g^0(\lambda))^{-1}
\]
determine a solution of the constraint optimisation problem
\begin{multline}
   \max_{(f,g)\in \mathcal{D}_f\times\md D_g}\bigg(\bigg\langle [D_N^{\mu}\me a_N]_{+\mu n}- \me T^{\mu}_N\me
    a_{\mu},
        (\me P^{\mu}_N)^{-1}
   [D_N^{\mu}\me a_N]_{+\mu n}
    -
   (\me P^{\mu}_N)^{-1}\me T^{\mu}_N\me a^{\mu}_N\bigg\rangle
       +\langle\me Q_N\me
    a_N,\me a_N\rangle\bigg)\\
    =\bigg\langle [D_N^{\mu}\me a_N]_{+\mu n}- (\me T^{\mu}_N)^0\me
    a^{\mu}_N,((\me P^{\mu}_N)^0)^{-1}[D_N^{\mu}\me a_N]_{+\mu n}
   - ((\me P^{\mu}_N)^0)^{-1}(\me T^{\mu}_N)^0\me a^{\mu}_N
    \bigg\rangle
    +\langle\me Q^0_N\me
    a_N,\me a_N\rangle.
\label{minimax1}
\end{multline}
The minimax spectral characteristic $h^0=h_{\mu,N}(f^0,g^0)$ is calculated by formula (18) if
$h_{\mu,N}(f^0,g^0)\in H_{\mathcal{D}}$.
\end{lema}

For more detailed analysis of properties of the least favorable spectral densities and minimax-robust spectral characteristics we observe that
the minimax spectral characteristic $h^0$ and the least favourable spectral densities $(f^0,g^0)$
form a saddle
point of the function $\Delta(h;f,g)$ on the set
$H_{\mathcal{D}}\times\mathcal{D}$.

The saddle point inequalities
\[\Delta(h;f^0,g^0)\geq\Delta(h^0;f^0,g^0)\geq\Delta(h^0;f,g)
\]
$ \forall f\in \mathcal{D}_f,\forall g\in \mathcal{D}_g,\forall h\in H_{\mathcal{D}}$ hold true if $h^0=h_{\mu,N}(f^0,g^0)$ and
$h_{\mu,N}(f^0,g^0)\in H_{\mathcal{D}}$, where $(f^0,g^0)$  is a
solution of the  constraint optimisation problem
\be  \label{cond-extr1}
\widetilde{\Delta}(f,g)=-\Delta(h_{\mu,N}(f^0,g^0);f,g)\to
\inf, (f,g)\in \mathcal{D},\ee
where the functional $\Delta(h_{\mu,N}(f^0,g^0);f,g)$ is calculated by the formula

\begin{align*}
\nonumber
&\Delta(h_{\mu,N}(f^0,g^0);f,g)=
\\
&\frac{1}{2\pi}\int_{-\pi}^{\pi}
\frac{\lambda^{2n}}{|1-e^{i\lambda\mu}|^{2n}}
\left[{(1-e^{i\lambda\mu})^n}(\vec{A}_N(e^{i\lambda}))^{\top}g^0(\lambda) +
(\vec{C}^0_{\mu,N}(e^{i \lambda}))^{\top}
\right](f^0(\lambda)+{\lambda}^{2n}g^0(\lambda))^{-1}
\times
\\
&\times
 f(\lambda)\, (f^0(\lambda)+{\lambda}^{2n}g^0(\lambda))^{-1}
\left[{(1-e^{-i\lambda\mu})^n}{\vec{A}_N(e^{-i\lambda})}g^0(\lambda) +
\vec{C}^0_{\mu,N}(e^{-i \lambda})
\right]
d\lambda+
\\
&\frac{1}{2\pi}\int_{-\pi}^{\pi}
\frac{1}{|1-e^{i\lambda\mu}|^{2n}}
\left[{(1-e^{i\lambda\mu})^n}(\vec{A}_N(e^{i\lambda}))^{\top}f^0(\lambda) -(\lambda)^{2n}
(\vec{C}^0_{\mu,N}(e^{i \lambda}))^{\top}
\right](f^0(\lambda)+{\lambda}^{2n}g^0(\lambda))^{-1}
\\
&\times
 g(\lambda)\, (f^0(\lambda)+{\lambda}^{2n}g^0(\lambda))^{-1}
\left[{(1-e^{-i\lambda\mu})^n}{\vec{A}_N(e^{-i\lambda})}f^0(\lambda) -(\lambda)^{2n}
\vec{C}^0_{\mu,N}(e^{-i \lambda})
\right]
d\lambda,
\end{align*}
where
\[
\vec{C}^0_{\mu,N}(e^{i \lambda})=\sum_{k=0}^{\infty}
(((\me P^{\mu}_N)^0)^{-1}[D_N^{\mu}\me a_N]_{+\mu n}-((\me P^{\mu}_N)^0)^{-1}(\me T^{\mu}_N)^0\me a^{\mu}_N)_k e^{ik\lambda}.
\]

The constrained optimisation problem (31) is equivalent to the unconstrained optimisation problem
\be  \label{uncond-extr}
\Delta_{\mathcal{D}}(f,g)=\widetilde{\Delta}(f,g)+ \delta(f,g|\mathcal{D}_f\times
\mathcal{D}_g)\to\inf,\ee
 where $\delta(f,g|\mathcal{D}_f\times
\mathcal{D}_g)$ is the indicator function of the set
$\mathcal{D}=\mathcal{D}_f\times\mathcal{D}_g$.
 Solution $(f^0,g^0)$ to this unconstrained optimisation problem is characterized by the condition $0\in
\partial\Delta_{\mathcal{D}}(f^0,g^0)$, where
$\partial\Delta_{\mathcal{D}}(f^0,g^0)$ is the subdifferential of the functional $\Delta_{\mathcal{D}}(f,g)$ at point $(f^0,g^0)\in \mathcal{D}=\mathcal{D}_f\times\mathcal{D}_g$.
 This condition makes it possible to find the least favourable spectral densities in some special classes of spectral densities $\mathcal{D}=\mathcal{D}_f\times\mathcal{D}_g$.

The form of the functional $\Delta(h_{\mu,N}(f^0,g^0);f,g)$ is convenient for application the Lagrange method of indefinite multipliers for
finding solution to the problem (32).
Making use of the method of Lagrange multipliers and the form of
subdifferentials of the indicator functions $\delta(f,g|\mathcal{D}_f\times
\mathcal{D}_g)$ of the set
$\mathcal{D}_f\times\mathcal{D}_g$ of spectral densities
we describe relations that determine least favourable spectral densities in some special classes
of spectral densities (see [20, 23] for additional details).

\subsection{Least favorable spectral density in classes $\md D_0 \times \md D_{1\delta}$}\label{set1}

Consider the problem of optimal linear estimation of the functional $A_N\vec{\xi}$
 which depends on unobserved values of a sequence $\vec\xi(m)$ with stationary increments based on observations of the sequence $\vec\xi(m)+\vec\eta(m)$
 at points  of the set $Z\setminus\{0,1,2,\ldots,N\}$
  under the condition that the sets of admissible spectral densities $\md D_{f0}^k, \md D_{1\delta}^{k},k=1,2,3,4$ are defined as follows:
$$\md D_{f0}^{1} =\left\{f\left|\frac{1}{2\pi} \int
_{-\pi}^{\pi}
\frac{|1-e^{i\lambda\mu}|^{2n}}{|\lambda|^{2n}}
f(\lambda )d\lambda  =P\right.\right\},$$
$$\md D_{f0}^{2} =\left\{f\left|\frac{1}{2\pi }
\int _{-\pi }^{\pi}
\frac{|1-e^{i\lambda\mu}|^{2n}}{|\lambda|^{2n}}
{\rm{Tr}}\,[ f(\lambda )]d\lambda =p\right.\right\},$$
$$\md D_{f0}^{3} =\left\{f\left|\frac{1}{2\pi }
\int _{-\pi}^{\pi}
\frac{|1-e^{i\lambda\mu}|^{2n}}{|\lambda|^{2n}}
f_{kk} (\lambda )d\lambda =p_{k}\right.\right\},$$
$$\md D_{f0}^{4} =\left\{f\left|\frac{1}{2\pi} \int _{-\pi}^{\pi}
\frac{|1-e^{i\lambda\mu}|^{2n}}{|\lambda|^{2n}}
\left\langle B_{1} ,f(\lambda )\right\rangle d\lambda  =p\right.\right\},$$
\begin{equation*}
\md D_{1\delta}^{1}=\left\{g\biggl|\frac{1}{2\pi} \int_{-\pi}^{\pi}
\left|{\rm{Tr}}(g(\lambda )-g_{1} (\lambda))\right|d\lambda \le \delta\right\},
\end{equation*}
\begin{equation*}
\md D_{1\delta}^{2}=\left\{g\biggl|\frac{1}{2\pi } \int_{-\pi}^{\pi}
\left|g_{kk} (\lambda )-g_{kk}^{1} (\lambda)\right|d\lambda  \le \delta_{k}\right\},
\end{equation*}
\begin{equation*}
\md D_{1\delta}^{3}=\left\{g\biggl|\frac{1}{2\pi } \int_{-\pi}^{\pi}
\left|\left\langle B_{2} ,g(\lambda )-g_{1}(\lambda )\right\rangle \right|d\lambda  \le \delta\right\},
\end{equation*}
\begin{equation*}
\md D_{1\delta}^{4}=\left\{g\biggl|\frac{1}{2\pi} \int_{-\pi}^{\pi}
\left|g_{ij} (\lambda )-g_{ij}^{1} (\lambda)\right|d\lambda  \le \delta_{i}^j\right\}.
\end{equation*}

\noindent
Here  $g_{1}(\lambda )$ is a fixed spectral density, $p, p_k, k=\overline{1,T}$, $\delta,\delta_{k},k=\overline{1,T}$, $\delta_{i}^{j}, i,j=\overline{1,T}$, are given numbers, $P, B_1, B_2$ are given positive-definite Hermitian matrices.

From the condition $0\in\partial\Delta_{\mathcal{D}}(f^0,g^0)$
we find the following equations which determine the least favourable spectral densities for these given sets of admissible spectral densities.

For the first set of admissible spectral densities $\md D_{f0}^1 \times \md D_{1\delta}^{1}$ we have equations

\begin{multline} \label{eq_4_1f}
\left(
{(1-e^{i\lambda\mu})^n}\vec{A}_N(e^{i\lambda})g^0(\lambda) +
\vec{C}^0_{\mu,N}(e^{i \lambda})
\right)
\left(
{(1-e^{i\lambda\mu})^n}\vec{A}_N(e^{i\lambda})g^0(\lambda) +
\vec{C}^0_{\mu,N}(e^{i \lambda})
\right)^{*}=
\\
=\left(\frac{|1-e^{i\lambda\mu}|^{2n}}{|\lambda|^{2n}} (f^0(\lambda)+{\lambda}^{2n}g^0(\lambda))\right)
\vec{\alpha}_f\cdot \vec{\alpha}_f^{*}\left(\frac{|1-e^{i\lambda\mu}|^{2n}}{|\lambda|^{2n}} (f^0(\lambda)+{\lambda}^{2n}g^0(\lambda))\right),
\end{multline}
\begin{multline} \label{eq_4_1g}
\left(
{(1-e^{i\lambda\mu})^n}\vec{A}_N(e^{i\lambda})f^0(\lambda) -(\lambda)^{2n}
\vec{C}^0_{\mu,N}(e^{i \lambda})
\right)
\times
 \\
\times
\left(
{(1-e^{i\lambda\mu})^n}\vec{A}_N(e^{i\lambda})f^0(\lambda) -(\lambda)^{2n}
\vec{C}^0_{\mu,N}(e^{i \lambda})
\right)^{*}=
\\
=
\beta^{2} \gamma_2( \lambda ){|1-e^{i\lambda\mu}|^{2n}}\left( f^0(\lambda)+{\lambda}^{2n}g^0(\lambda)\right)^2,
\end{multline}
\begin{equation} \label{eq_4_1c}
\frac{1}{2 \pi} \int_{-\pi}^{ \pi}
\left|{\mathrm{Tr}}\, (g_0( \lambda )-g_{1}(\lambda )) \right|d\lambda =\delta,
\end{equation}
\noindent where $\vec{\alpha}_f$, $\beta^{2}$,  are  Lagrange multipliers, the function $\left| \gamma_2( \lambda ) \right| \le 1$ and
\[\gamma_2( \lambda )={ \mathrm{sign}}\; ({\mathrm{Tr}}\, (g_{0} ( \lambda )-g_{1} ( \lambda ))): \; {\mathrm{Tr}}\, (g_{0} ( \lambda )-g_{1} ( \lambda )) \ne 0.\]

For the second set of admissible spectral densities $\md D_{f0}^2 \times \md D_{1\delta}^{2}$ we have equations
\begin{multline} \label{eq_4_2f}
\left(
{(1-e^{i\lambda\mu})^n}\vec{A}_N(e^{i\lambda})g^0(\lambda) +
\vec{C}^0_{\mu,N}(e^{i \lambda})
\right)
\left(
{(1-e^{i\lambda\mu})^n}\vec{A}_N(e^{i\lambda})g^0(\lambda) +
\vec{C}^0_{\mu,N}(e^{i \lambda})
\right)^{*}=
\\
=\alpha_f^{2} \left(\frac{|1-e^{i\lambda\mu}|^{2n}}{|\lambda|^{2n}} (f^0(\lambda)+{\lambda}^{2n}g^0(\lambda))\right)^2,
\end{multline}
\begin{multline}   \label{eq_4_2g}
\left(
{(1-e^{i\lambda\mu})^n}\vec{A}_N(e^{i\lambda})f^0(\lambda) -(\lambda)^{2n}
\vec{C}^0_{\mu,N}(e^{i \lambda})
\right)
\times
 \\
\times
\left(
{(1-e^{i\lambda\mu})^n}\vec{A}_N(e^{i\lambda})f^0(\lambda) -(\lambda)^{2n}
\vec{C}^0_{\mu,N}(e^{i \lambda})
\right)^{*}=
\\
\left(|1-e^{i\lambda\mu}|^{n}(f^0(\lambda)+{\lambda}^{2n}g^0(\lambda))\right)
\left \{ \beta_{k}^{2} \gamma^2_{k} ( \lambda ) \delta_{kl} \right \}_{k,l=1}^{T}
\left(|1-e^{i\lambda\mu}|^{n}(f^0(\lambda)+{\lambda}^{2n}g^0(\lambda))\right),
\end{multline}
\begin{equation} \label{eq_4_2c}
\frac{1}{2 \pi} \int_{- \pi}^{ \pi}  \left|g^0_{kk} ( \lambda)-g_{kk}^{1} ( \lambda ) \right| d\lambda =\delta_{k},
\end{equation}

\noindent where $\alpha _{f}^{2}$, $\beta_{k}^{2}$ are Lagrange multipliers,  the
functions $\left| \gamma^2_{k} ( \lambda ) \right| \le 1$ and
\[\gamma_{k}^2( \lambda )={ \mathrm{sign}}\;(g_{kk}^{0}( \lambda)-g_{kk}^{1} ( \lambda )): \; g_{kk}^{0} ( \lambda )-g_{kk}^{1}(\lambda ) \ne 0, \; k= \overline{1,T}.\]

For the third set of admissible spectral densities $\md D_{f0}^{3}\times \md D_{1\delta}^{3}$ we have equations
\begin{multline} \label{eq_4_3f}
\left(
{(1-e^{i\lambda\mu})^n}\vec{A}_N(e^{i\lambda})g^0(\lambda) +
\vec{C}^0_{\mu,N}(e^{i \lambda})
\right)
\left(
{(1-e^{i\lambda\mu})^n}\vec{A}_N(e^{i\lambda})g^0(\lambda) +
\vec{C}^0_{\mu,N}(e^{i \lambda})
\right)^{*}=
\\
\left(\frac{|1-e^{i\lambda\mu}|^{2n}}{|\lambda|^{2n}} (f^0(\lambda)+{\lambda}^{2n}g^0(\lambda))\right)
\left\{\alpha _{fk}^{2} \delta _{kl} \right\}_{k,l=1}^{T}
\left(\frac{|1-e^{i\lambda\mu}|^{2n}}{|\lambda|^{2n}} (f^0(\lambda)+{\lambda}^{2n}g^0(\lambda))\right),
\end{multline}
\begin{multline}   \label{eq_4_3g}
\left(
{(1-e^{i\lambda\mu})^n}\vec{A}_N(e^{i\lambda})f^0(\lambda) -(\lambda)^{2n}
\vec{C}^0_{\mu,N}(e^{i \lambda})
\right)
\times
 \\
\times
\left(
{(1-e^{i\lambda\mu})^n}\vec{A}_N(e^{i\lambda})f^0(\lambda) -(\lambda)^{2n}
\vec{C}^0_{\mu,N}(e^{i \lambda})
\right)^{*}=
\\
=
\beta^{2} \gamma_2'( \lambda )
\left(|1-e^{i\lambda\mu}|^{n}(f^0(\lambda)+{\lambda}^{2n}g^0(\lambda))\right)
B_{2}^{ \top}
\left(|1-e^{i\lambda\mu}|^{n}(f^0(\lambda)+{\lambda}^{2n}g^0(\lambda))\right),
\end{multline}
\begin{equation} \label{eq_4_3c}
\frac{1}{2 \pi} \int_{- \pi}^{ \pi} \left| \left \langle B_{2}, g_0( \lambda )-g_{1} ( \lambda ) \right \rangle \right|d\lambda
= \delta,
\end{equation}
where $\alpha_{fk}^{2}$,  $\beta^{2}$ are Lagrange multipliers, the function
$\left| \gamma_2' ( \lambda ) \right| \le 1$ and
\[\gamma_2' ( \lambda )={ \mathrm{sign}}\; \left \langle B_{2},g_{0} ( \lambda )-g_{1} ( \lambda ) \right \rangle : \; \left \langle B_{2},g_{0} ( \lambda )-g_{1} ( \lambda ) \right \rangle \ne 0.\]

For the fourth set of admissible spectral densities $\md D_{f0}^4 \times\md D_{1\delta}^{4}$ we have equations
\begin{multline} \label{eq_4_4f}
\left(
{(1-e^{i\lambda\mu})^n}\vec{A}_N(e^{i\lambda})g^0(\lambda) +
\vec{C}^0_{\mu,N}(e^{i \lambda})
\right)
\left(
{(1-e^{i\lambda\mu})^n}\vec{A}_N(e^{i\lambda})g^0(\lambda) +
\vec{C}^0_{\mu,N}(e^{i \lambda})
\right)^{*}=
\\
=
\alpha_f^{2} \left(\frac{|1-e^{i\lambda\mu}|^{2n}}{|\lambda|^{2n}} (f^0(\lambda)+{\lambda}^{2n}g^0(\lambda))\right)
B_{1}^{\top}
\left(\frac{|1-e^{i\lambda\mu}|^{2n}}{|\lambda|^{2n}} (f^0(\lambda)+{\lambda}^{2n}g^0(\lambda))\right),
\end{multline}
\begin{multline}  \label{eq_4_4g}
\left(
{(1-e^{i\lambda\mu})^n}\vec{A}_N(e^{i\lambda})f^0(\lambda) -(\lambda)^{2n}
\vec{C}^0_{\mu,N}(e^{i \lambda})
\right)
\times
 \\
\times
\left(
{(1-e^{i\lambda\mu})^n}\vec{A}_N(e^{i\lambda})f^0(\lambda) -(\lambda)^{2n}
\vec{C}^0_{\mu,N}(e^{i \lambda})
\right)^{*}=
\\
=
\left(|1-e^{i\lambda\mu}|^{n}(f^0(\lambda)+{\lambda}^{2n}g^0(\lambda))\right)
\left \{ \beta_{ij}( \lambda ) \gamma_{ij} ( \lambda ) \right \}_{i,j=1}^{T}
\left(|1-e^{i\lambda\mu}|^{n}(f^0(\lambda)+{\lambda}^{2n}g^0(\lambda))\right),
\end{multline}
\begin{equation} \label{eq_4_4c}
\frac{1}{2 \pi} \int_{- \pi}^{ \pi}  \left|g^0_{ij}(\lambda)-g_{ij}^{1}( \lambda ) \right|d\lambda = \delta_{i}^{j},
\end{equation}
where $\alpha _{f}^{2}$, $ \beta_{ij}$ are Lagrange multipliers, the functions $\left| \gamma_{ij} ( \lambda ) \right| \le 1$ and
\[
\gamma_{ij} ( \lambda )= \frac{g_{ij}^{0} ( \lambda )-g_{ij}^{1} (\lambda )}{ \left|g_{ij}^{0} ( \lambda )-g_{ij}^{1}(\lambda) \right|}: \; g_{ij}^{0} ( \lambda )-g_{ij}^{1} ( \lambda ) \ne 0, \; i,j= \overline{1,T}.
\]

The following theorem  holds true.

\begin{thm}
Let the minimality condition (12) hold true. The least favorable spectral densities $f_{0}(\lambda)$, $g_{0}(\lambda)$ in classes
$\md D_{f0}^{k}\times \md D_{1\delta}^{k},k=1,2,3,4$ for the optimal linear estimation of the functional  $A_N\vec{\xi}$ from observations of the sequence $\vec{\xi}(m)+ \vec{\eta}(m)$ at points  of the set $Z\setminus\{0,1,2,\ldots,N\}$
are determined by equations
(33) - (35), (36) - (38), (39) - (41), (42) - (44),
respectively,
the constrained optimization problem (30) and restrictions  on densities from the corresponding classes
$ \md D_{f0}^{k}, \md D_{1\delta}^{k},k=1,2,3,4$.  The minimax-robust spectral characteristic of the optimal estimate of the functional $A_N\vec{\xi}$ is determined by the formula (18).
\end{thm}

\section{Conclusions}

In this article, we present results of investigation of stochastic sequences with periodically stationary increments.
We give definition of increment sequence and introduce stochastic sequences  with periodically stationary (periodically correlated, cyclostationary) increments.
These non-stationary stochastic sequences combine  periodic structure of covariation functions of sequences as well as integrating one.

We describe methods of solution of the classical estimation problem for linear functionals  constructed from unobserved values of a sequence with periodically stationary increments.
Estimates are based on observations of the sequence with a stationary noise sequence.
Estimates are obtained by representing the sequence under investigation as a vector-valued sequence with stationary increments.
The problem is investigated in the case of spectral certainty, where spectral densities of sequences are exactly known.
In this case  we propose an approach based on the Hilbert space projection method.
We derive formulas for calculating the spectral characteristics and the mean-square errors of the optimal estimates of the functionals.
In the case of spectral uncertainty where the spectral densities are not exactly known while, instead, some sets of admissible spectral densities are specified,
the minimax-robust method is applied.
We propose a representation of the mean square error in the form of a linear
functional in $L_1$ space with respect to spectral densities, which allows
us to solve the corresponding constrained optimization problem and
describe the minimax-robust estimates of the functionals. Formulas
that determine the least favorable spectral densities and minimax-robust spectral characteristic of the optimal linear estimates of
the functionals are derived  for a collection of specific classes
of admissible spectral densities.

These least favourable spectral density matrices are solutions of the optimization problem $\Delta_{\mathcal{D}}(f,g)=\widetilde{\Delta}(f,g)+ \delta(f,g|\mathcal{D}_f\times
\mathcal{D}_g)\to\inf,$
 where $\delta(f,g|\mathcal{D}_f\times
\mathcal{D}_g)$ is the indicator function of the set
$\mathcal{D}=\mathcal{D}_f\times\mathcal{D}_g$.
 Solution $(f^0,g^0)$ to this unconstraint optimisation problem is characterized by the condition $0\in
\partial\Delta_{\mathcal{D}}(f^0,g^0)$, where
$\partial\Delta_{\mathcal{D}}(f^0,g^0)$ is the subdifferential of the functional $\Delta_{\mathcal{D}}(f,g)$ at point $(f^0,g^0)\in \mathcal{D}=\mathcal{D}_f\times\mathcal{D}_g$.
 This condition makes it possible to find the least favourable spectral densities in some special classes of spectral densities.
These are: classes $D_0$ of densities with the moment restrictions,
 classes $D_{1\delta}$ which describe the ``$\delta$-neighborhood''\, models in the space $L_1$ of a fixed bounded spectral density.

\def\selectlanguageifdefined#1{
\expandafter\ifx\csname date#1\endcsname\relax
\else\selectlanguage{#1}\fi}
\providecommand*{\href}[2]{{\small #2}}
\providecommand*{\url}[1]{{\small #1}}
\providecommand*{\BibUrl}[1]{\url{#1}}
\providecommand{\BibAnnote}[1]{}
\providecommand*{\BibEmph}[1]{#1}

\renewcommand{\refname}{References}

\end{document}